\documentclass[12pt,a4paper]{article}%
\usepackage{amsmath}
\usepackage{graphicx}
\usepackage{epsfig}%
\usepackage{amsfonts}%
\usepackage{amssymb}
\usepackage{enumitem}
\usepackage{float}
\usepackage{tikz}
\usepackage{cleveref}
\usepackage{color}

\newcommand{\cred}{\color{red}}
\newcommand{\cn}{\color{black}}
    \usepackage[numbers,sort]{natbib}
\usepackage[normalem]{ulem}
 \usepackage{authblk}
 
\newtheorem{theorem}{Theorem}[section]

\newtheorem{proposition}[theorem]{Proposition}

\makeatletter\@addtoreset{equation}{section}\makeatother
\makeatletter\@addtoreset{figure}{section}\makeatother
\makeatletter\@addtoreset{table}{section}\makeatother
\newenvironment{rcases}
  {\left.\begin{aligned}}
  {\end{aligned}\right\rbrace}
\RequirePackage[%
    font={small,sf},
    labelfont=bf,
    format=hang,    
    format=plain,
    margin=0pt,
    width=0.8\textwidth,
]{caption}

\begin{document}


\title{Computing electrostatic potentials using regularization based on the range-separated 
tensor format}
 
\author[1]{Peter Benner \thanks{\tt benner@mpi-magdeburg.mpg.de}}
\author[1,2]{Venera Khoromskaia \thanks{\tt vekh@mis.mpg.de}}
\author[2]{Boris Khoromskij \thanks{\tt bokh@mis.mpg.de}}
\author[1]{Cleophas Kweyu \thanks{\tt kweyu@mpi-magdeburg.mpg.de}}
\author[1]{Matthias Stein \thanks{\tt matthias.stein@mpi-magdeburg.mpg.de}}
\affil[1]{Max Planck Institute for Dynamics of Complex Technical Systems, Sandtorstr.~1, D-39106 Magdeburg, 
Germany}
\affil[2]{Max Planck Institute for Mathematics in the Sciences, Inselstr.~22-26, D-04103 Leipzig, Germany}     

\date{\vspace{-7ex}}

\maketitle

\begin{abstract}
In this paper, we apply the range-separated (RS) tensor format \cite{BKK_RS:18} for the construction 
of new regularization scheme  for the Poisson-Boltzmann equation (PBE) describing the electrostatic 
potential in biomolecules. In our approach, we use the RS tensor representation to   the  discretized Dirac 
delta \cite{khor-DiracRS:2018} to construct  an   efficient RS splitting of the PBE solution in the solute 
(molecular) region. The PBE then needs to be solved with a regularized source term, and thus black-box 
solvers can be applied. The main computational benefits are due to the localization of the modified 
right-hand side within the molecular region and automatic maintaining of the continuity in the Cauchy data 
on the interface. Moreover, this computational scheme only includes solving a single system of FDM/FEM 
equations for the smooth long-range (i.e., regularized) part of the collective potential represented by a 
low-rank RS-tensor with a controllable precision. The total potential is obtained by adding this solution 
to the directly precomputed rank-structured tensor representation for the short-range contribution. 
Enabling finer grids in PBE computations is another advantage of the proposed techniques. In the numerical 
experiments, we consider only the free space electrostatic potential for proof of concept. We illustrate 
that the classical Poisson equation (PE) model does not accurately capture the solution singularities in 
the numerical approximation as compared to the new approach by the RS tensor format.
\end{abstract}



\section{Introduction}
\label{sec:Intro}

Numerical treatment of long-range interaction potentials is a challenging task in computer modeling of 
multiparticle systems, for example,  in calculation of electrostatics in large solvated biological 
systems, in protein docking, or in many particles dynamics simulations
\cite{PolGlos:96,DesHolm:98,HuMcCam:1999,Stein:2010,LiStCaMaMe:13}.
The well-known Poisson-Boltzmann equation (PBE) introduced and analyzed for example, in \cite{Holst94}, 
is one of the most popular implicit solvent models for computation of the electrostatic potential in 
proteins. Other models include the generalized Born (GB) methods \cite{Bashford:2000} and the 
polarizable continuum models (PCM) \cite{Barone:1997}. The PBE computes the electrostatic potential both 
in the protein and in the surrounding solvent, and it is widely used in protein docking, in 
classification problems, and for computation of the free energy of biomolecules in a self-consistent way.

The main difficulty in the traditional finite element method (FEM) approximation schemes for the PBE problem
is related to the presence of a highly singular source term that includes a large sum of Dirac delta 
distributions which need to be resolved using rather coarse grids. To overcome these limitations, a number 
of regularization schemes for the finite element method (FEM) applied to the PBE, based on the full grid 
representation of all functional data, have been considered in the literature, see for example 
\cite{Holst2000,Xie:14} and references therein. Consequently, we note that the PBE theory has recently 
received major improvement in terms of accuracy by the introduction of solution decomposition techniques 
which have been developed for example, in \cite{Xie:14,Mirzadeh:13,Chen:07}, where the PBE is treated as 
an interface problem. This aims at avoiding the discontinuities in traces and fluxes at the interface 
between the biomolecule and the solvent and also to circumvent building the numerical approximations 
corresponding to the Dirac delta distributions because of the existence of an analytical expansions in the 
solute sub-region. 

However, these techniques still face the following computational challenges. First, jumps in the interface 
conditions, arising due to regularization splitting of the solution, need to be incorporated to eliminate 
the solution discontinuity (e.g., Cauchy data) at the interface. Second, the boundary conditions have to 
be specified using some analytical representation of the solution of the PBE. And third, in 
regularization-based techniques, see for example, \cite{Xie:14}, one has to solve multiple algebraic 
systems for the linear and nonlinear boundary value problems before summing up the partial solutions which 
increases the computational costs. We provide an overview of these techniques in  \Cref{Soln_decomp}.

Here, we present a new approach for the regularization of the PBE by using the range-separated (RS) 
canonical tensor format introduced and analyzed in \cite{BKK_RS:18}. The RS tensor format relies on the 
independent grid-based low-rank tensor representation of the long- and short-range parts in the total sum 
of single-particle electrostatic potentials discretized on a fine 3D $n\times n\times n$ Cartesian grid 
$\Omega_n$ in the computational box $\Omega \subset \mathbb{R}^3$. This representation is based on the 
splitting of a single reference potential, defined by a radial function like 
$p(\|\bar{x}\|)=1/\|\bar{x}\|$, into a sum of localized and long-range low-rank canonical tensors 
both represented on the computational grid $\Omega_n$. The long-range part in the collective 
potential of a many-particle system is represented as a low-rank canonical tensor with the rank only 
logarithmically depending on the number of particles in the system. The rank reduction algorithm is 
performed by the canonical-to-Tucker (C2T) transform via the reduced higher order singular value 
decomposition (RHOSVD) \cite{khor-ml-2009} with a subsequent Tucker-to-canonical (T2C) decomposition. 
The short-range contributions to a many-particle potential are parametrized by a single low-rank 
canonical tensor of local support. Notice that in \cite{BKK_RS:16} \cred,\cn it was already 
sketched how the RS tensor 
formats may be utilized for calculation of the free interaction energy of 
protein-type systems, and the idea for regularized formulation of the PBE by 
using the smooth long-range part of the free space electrostatic potential was 
outlined. 

In this paper \cred,\cn we introduce the new regularization scheme for the solution of the PBE
adapting the RS tensor format. It is based on a certain splitting scheme for 
the highly singular solution and right-hand side in PBE, by using 
the RS tensor decomposition of the discretized Dirac delta introduced 
in \cite{khor-DiracRS:2018}. 
This approach requires only a simple modification (regularization) in 
the right-hand side of the PBE in the solute region, but it does not change the interface conditions and, 
hence, the FEM system matrix. The most singular component in the potential is recovered explicitly by the 
short range part in the RS tensor splitting of the free space potential. The main computational benefits 
are due to the localization of the modified right-hand side within the molecular region and automatic 
maintaining of the continuity in the Cauchy data on the interface. Furthermore, this computational scheme 
only includes solving a single system of FEM equations for the smooth long-range (i.e., regularized) part 
of the collective potential represented by a low-rank RS-tensor with a controllable precision. The total 
potential is obtained by adding this solution to the directly precomputed rank-structured tensor 
representation for the short-range contribution. 
 
As numerical illustrations, we compute the free-space electrostatic potentials of biomolecules using the 
RS  tensor format in the framework of splitting scheme, and compare them with the solutions calculated by 
the traditional FEM/FDM discretization methods for the PBE \cite{Holst94,Holst2000,Bakersept2001}. 

The rest of the paper is organized as follows. Section \ref{Soln_decomp} provides a short overview of 
existing solution decomposition schemes for the PBE problem. Section \ref{sec:RS_survey} describes the 
principles of the rank-structured tensor approximation to the long-range electrostatic potential and 
sketches the RS tensor decomposition techniques \cite{BKK_RS:18} for the free space electrostatic potential 
of many particle systems. The main Section \ref{sec:RS_2_PBE} explains how the application of the RS tensor 
format leads to the new regularization scheme  for solving the PBE. 
Finally, Section \ref{sec:Numer_Tests} 
presents the numerical tests illustrating the benefits of the proposed method, and comparisons with the 
solutions obtained by the standard FDM-based PBE solver are provided.

\section{On existing solution decomposition techniques}\label{Soln_decomp}
\begin{figure}[b]
\centering
\includegraphics[height=6.0cm]{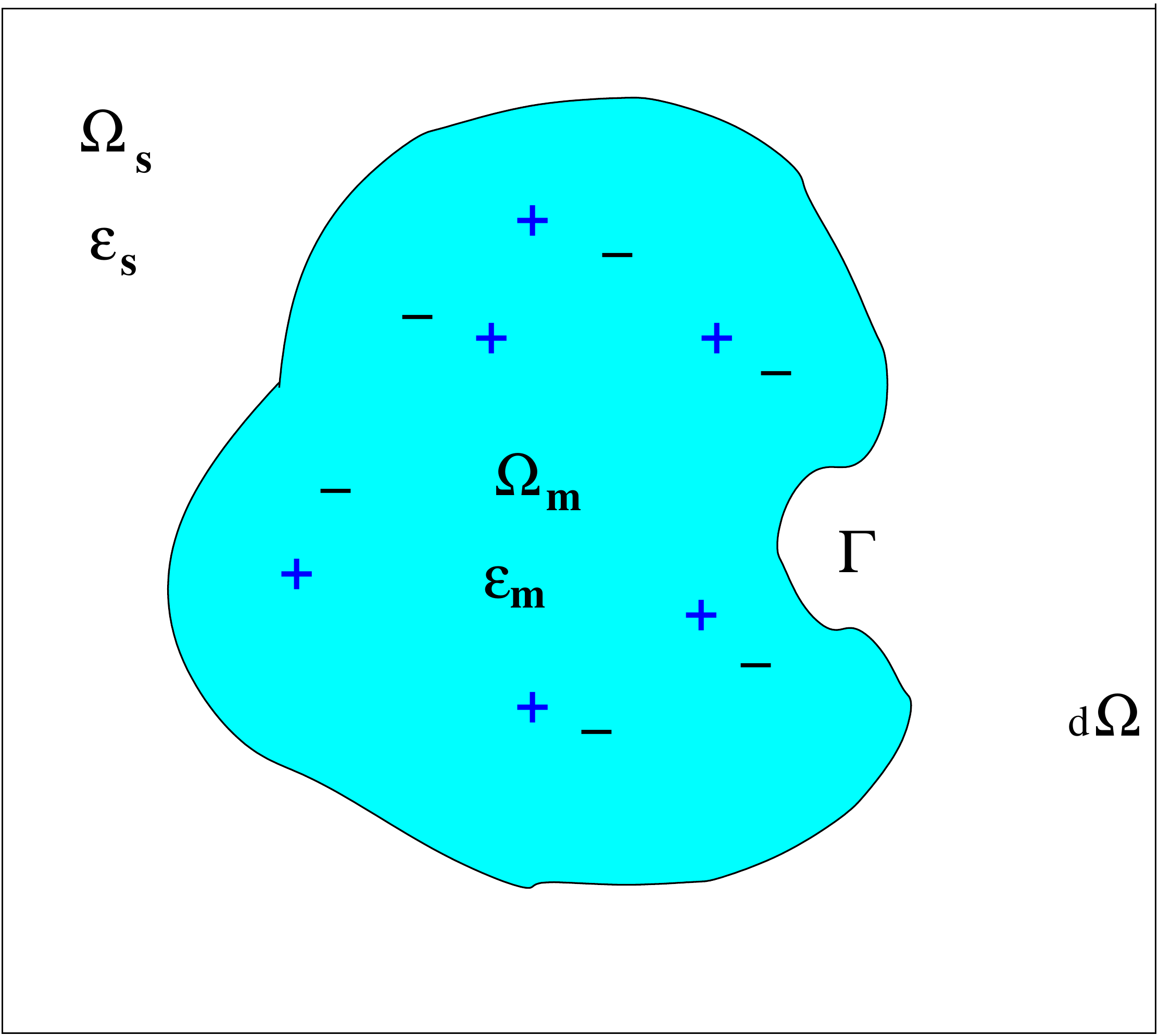} 
\caption{Solute and solvent regions in the computational domain for the PBE. }
\label{fig:Protein}
\end{figure}

The PBE computes the dimensionless potential $u(\bar{x}) = {e_c\psi(\bar{x})}/{K_B T}$, which is  
scaled by $e_c/{K_B T}$, and $\psi(\bar{x})$ is the original electrostatic potential in 
centimeter-gram-second (cgs) units at $ \bar{x}=(x,y,z) \in \mathbb{R}^3 $. It is given by
\begin{equation}\label{eq:PBE}
-\nabla\cdot(\epsilon(\bar{x})\nabla u(\bar{x})) + \bar{\kappa}^2(\bar{x})\sinh(u(\bar{x})) = 
  \frac{4\pi e_c^2}{K_B T}\sum_{i=1}^{N_m}z_i\delta(\bar{x}-\bar{x}_i),\quad \Omega \in 
\mathbb{R}^3,
\end{equation}
subject to 
\begin{equation}\label{eq:DH_solution}
u(\bar{x}) = \frac{e_c^2}{K_B T} \sum_{i=1}^{N_m}\frac{z_ie^{-\bar{\kappa}(d-a_i)}}{\epsilon_s 
(1+\bar{\kappa} a_i)d} \quad \mbox{on} \: \partial{\Omega}, \quad d = \lVert \bar{x}-\bar{x}_i \rVert,      
\end{equation}
where $K_B T$, $K_B$, $T$, and $e_c$ are the thermal energy, the Boltzmann constant, the absolute 
temperature, and the electron charge, respectively, $\bar{\kappa}^2 = {8\pi e_c^2 I}/{1000\epsilon_s K_BT}$ 
is a function of the ionic strength $I = 1/2\sum_{j=1}^{m}c_jz_j^2$, where $c_j$ and $z_j$ are the charge 
and concentration of each ion. The sum of Dirac delta distributions, located at atomic centers $\bar{x}_i$, 
represents the molecular charge density, $z_i$ are the point partial charges of the protein, $\epsilon_s$ 
is the solvent dielectric constant, $a_i$ are the atomic radii, and $N_m$ is the total number of point partial 
charges in the protein. The functions $\epsilon(\bar{x})$ and $\bar{\kappa}^2(\bar{x})$ are piecewise 
constant defined by
\begin{eqnarray}\label{eqn:diel_kappa_def}
 \epsilon(\bar{x}) =
  \begin{cases}
   \epsilon_m = 2 & \text{if } \bar{x} \in \Omega_m\\
   \epsilon_s \,\,= 78.54 & \text{if } \bar{x} \in \Omega_s
  \end{cases}, \quad \quad
 \bar{\kappa}(\bar{x}) =
  \begin{cases}
   0 & \text{if } \bar{x} \in \Omega_m \\
    \sqrt{\epsilon_s}\bar{\kappa} & \text{if } \bar{x} \in \Omega_s 
  \end{cases},
\end{eqnarray}
where $\Omega_m$ and $\Omega_s$ are the molecular and solvent regions, respectively, as shown in 
\Cref{fig:Protein}.

In order to overcome the difficulties arising from the singularities caused by the impulsive source term,
several solution decomposition techniques have been suggested in the literature. In the following, we will 
discuss those techniques that form the state of the art for solving the PBE. Following \cite{Xie:14}, the 
first solution decomposition is generally given by 
\begin{equation}\label{eq:PBE_interfaceI}
\begin{rcases}
\begin{aligned}
 -\epsilon_m \Delta u(\bar{x}) &= C \sum_{i=1}^{N_m}z_i\delta(\bar{x}-\bar{x}_i),  \quad \quad  
				  &\bar{x} \in \Omega_m,\\
 -\epsilon_s \Delta u(\bar{x}) + \bar{\kappa}^2 \sinh(u(\bar{x})) &= 0,  & \bar{x} \in \Omega_s,\\
 u(s^+) = u(s^-), \quad \epsilon_s \frac{\partial u(s^+)}{\partial n(s)} 
 &= \epsilon_m \frac{\partial u(s^-)}{\partial n(s)}, & s \in \varGamma,\\
 u(s) &= g(s), & s \in \partial \Omega,
\end{aligned}
\end{rcases}
\end{equation}
where $C = \frac{4\pi e_c^2}{K_B T}$, $s^+$ and $s^-$ represent the grid points in the vicinity of the 
interface $\varGamma$ in the solvent and the molecular regions, respectively. The solution $u(\bar{x})$ is 
decomposed as follows,
\begin{equation}
   u(\bar{x}) = G(\bar{x}) + \tilde{\phi}(\bar{x}) + \tilde{\psi}(\bar{x}).
 \end{equation}
The corresponding components of $u(\bar{x})$ include the analytical solution $G(\bar{x})$ of the Poisson 
equation in the molecular domain,
\begin{equation}\label{eq:PBE_solndecompIa}
    G(\bar{x}) = \frac{\hat{C}}{\epsilon_m}\sum_{i=1}^{N_m}\frac{z_i}{\| \bar{x}-\bar{x}_i \|}, \quad 
   \quad \bar{x} \in \Omega_m,
\end{equation}
the solution of the linear interface boundary value problem
\begin{equation}\label{eq:PBE_solndecompIb}
\begin{rcases}
\begin{aligned}
 \Delta \tilde{\phi}(\bar{x}) &= 0,  & \bar{x} \in \Omega_m \cup \Omega_s,\\
 \tilde{\phi}(s^+) &= \tilde{\phi}(s^-), \quad \epsilon_s \frac{\partial \tilde{\phi}(s^+)}{\partial n(s)} 
  = \epsilon_m \frac{\partial \tilde{\phi}(s^-)}{\partial n(s)} +  
 (\epsilon_m - \epsilon_s)\frac{G(s)}{\partial n(s)}, & s \in \varGamma,\\
 u(s) & = g(s) - G(s), & s \in \partial \Omega,
\end{aligned}
\end{rcases}
\end{equation}
and the solution of the nonlinear interface boundary value problem
\begin{equation}\label{eq:PBE_solndecompIc}
\begin{rcases}
\begin{aligned}
 \Delta \tilde{\psi}(\bar{x}) & = 0,  & \bar{x} \in \Omega_m,\\
 -\epsilon_s \Delta \tilde{\psi}(\bar{x}) + \bar{\kappa}^2 \sinh(\tilde{\psi}(\bar{x}) + \tilde{\phi}(\bar{x}) +
 G(\bar{x})) & = 0, & \bar{x} \in \Omega_s,\\
 \tilde{\psi}(s^+) = \tilde{\psi}(s^-), \quad \epsilon_s \frac{\partial \tilde{\psi}(s^+)}{\partial n(s)} 
 & = \epsilon_m \frac{\partial \tilde{\psi}(s^-)}{\partial n(s)}, & s \in \varGamma,\\
 u(s) & = 0, & s \in \partial \Omega.
\end{aligned}
\end{rcases}
\end{equation}
Secondly, we recall the solution decomposition from \cite{Mirzadeh:13} which takes the form 
\[u(\bar{x}) = \hat{u}(\bar{x}) + \tilde{u}(\bar{x}).\] The short-range part $\hat{u}(\bar{x})$ is given by 
\begin{eqnarray}\label{eq:PBE_solndecompIIa}
 \hat{u}(\bar{x}) =
  \begin{cases}
   G(\bar{x}) + u^0(\bar{x}) & \mbox{if} \quad \bar{x} \in \Omega_m,\\
   0 & \mbox{if} \quad \bar{x} \in \Omega_s,
  \end{cases}
\end{eqnarray}
where $u^0(\bar{x})$ is a harmonic function which compensates for the discontinuity on the interface and 
satisfies the following conditions
 \begin{equation}\label{eq:PBE_solndecompIIb}
\begin{rcases}
\begin{aligned}
 \Delta u^0(\bar{x}) & = 0,  & \mbox{if} \quad \bar{x} \in \Omega_m,\\
 u^0(s) &= -G(s), & s \in \varGamma.
\end{aligned}
\end{rcases}
\end{equation}
The regular part $\tilde{u}(\bar{x})$ is represented by
\begin{equation}\label{eq:PBE_solndecompIIc}
\begin{rcases}
\begin{aligned}
-\nabla \cdot(\epsilon(\bar{x})\nabla \tilde{u}(\bar{x})) + \bar{\kappa}^2(\bar{x})\sinh(\tilde{u}(\bar{x})) &= 0, \\
[\tilde{u}(\bar{x})]_{\varGamma} = 0, \quad [\epsilon \nabla \tilde{u}(\bar{x})\cdot \textbf{n}]_{\varGamma}
&= -\epsilon_m \nabla(G(\bar{x}) + u^0(\bar{x}))\cdot \textbf{n}\rvert_{\varGamma}.
\end{aligned}
\end{rcases}
\end{equation}
Lastly, the solution decomposition in \cite{Chen:07} is as follows: $u(\bar{x}) = G(\bar{x}) + u^r(\bar{x})$, where $G(\bar{x})$ 
is as in (\ref{eq:PBE_solndecompIa}) and the regular part is given by
\begin{equation}\label{eq:PBE_solndecompIIIa}
\begin{rcases}
\begin{aligned}
-\nabla \cdot(\epsilon\vec{\nabla}u^r) + \bar{\kappa}^2\sinh(u^r \!+\!G) 
   &= \nabla \cdot((\epsilon-\epsilon_m)\vec{\nabla}G),  & \, \textrm{in} \,\, \Omega \\
u^r &= g - G, & \, \textrm{on} \,\, \partial{\Omega}.
\end{aligned}
\end{rcases}
\end{equation}
We can also further decompose (\ref{eq:PBE_solndecompIIIa}) into the linear and nonlinear 
components so that $u^r(\bar{x}) = u^l(\bar{x}) + u^n(\bar{x})$, where
\begin{equation}\label{eq:PBE_solndecompIIIb}
\begin{rcases}
\begin{aligned}
-\nabla \cdot(\epsilon\vec{\nabla}u^l) &= \nabla \cdot((\epsilon-\epsilon_m)\nabla G),  
& \, \textrm{in} \,\, \Omega \\
u^l &= 0 & \, \textrm{on} \,\, \partial{\Omega},
\end{aligned}
\end{rcases}
\end{equation}
and
\begin{equation}\label{eq:PBE_solndecompIIIc}
\begin{rcases}
\begin{aligned}
-\nabla \cdot(\epsilon\vec{\nabla}u^n) + \bar{\kappa}^2\sinh(u^n + u^l + G) &= 0,  
& \, \textrm{in} \,\, \Omega \\
u^n &= g-G & \, \textrm{on} \,\, \partial{\Omega}.
\end{aligned}
\end{rcases}
\end{equation}

The fundamental idea in the above decomposition strategies is the pursuit of an efficient 
solution decomposition technique for the short- and long-range parts in a target tensor. However, 
all these techniques do not efficiently separate the long- and short-range components in the total 
potential sum, hence there is a need to apply so-called interface (or jump) conditions at the 
interface between the molecular region and the solvent region in order to reduce discontinuities.
In what follows, we describe the RS tensor format developed in \cite{BKK_RS:18}, applied as the main tool
in our paper to modify the PBE to increase the accuracy of its numerical approximation. A more detailed 
description of the new solution decomposition technique is provided in Sections \ref{ssec:PBE_Applic} 
and \ref{ssec:Disc_Split_scheme}.

\section{Rank-structured approximation of electrostatic potentials}
\label{sec:RS_survey}

Tensor-structured  numerical methods are now becoming popular in scientific computing due to their 
intrinsic property of reducing the grid-based solution of the multidimensional problems to 
essentially ``one-dimensional'' computations. These methods evolved from bridging of the traditional 
rank-structured tensor formats of multilinear algebra \cite{Kolda:07,smilde-book-2004} with the 
nonlinear approximation theory based on a separable representation of multidimensional functions and 
operators \cite{HaKhtens:04I,GHK:05,khor-rstruct-2006}. One of the ingredients in the development of 
tensor methods was the RHOSVD which allows to reduce the rank of tensors in a canonical format by the 
C2T decomposition without the need to construct the full size tensor \cite{khor-ml-2009}. Originally, 
it was used for the reduction of the ranks of canonical tensors when calculating three-dimensional 
convolution integrals in computational quantum chemistry, see the surveys 
\cite{VeKhorTromsoe:15,Khor-book-2018} and the references therein.

Recently, tensor-based approaches were suggested as new methods for the calculation of multiparticle 
long-range interaction potentials. For a given non-local generating kernel $p(\|\bar{x}\|)$, 
$\bar{x} \in \mathbb{R}^3$, the calculation of the weighted sum of interaction potentials in an 
$N$-particle system, with the particle locations at $\bar{x}_\nu \in \mathbb{R}^3$, $\nu=1,...,N$, 
\begin{equation}\label{eqn:PotSum1}
 P_0(\bar{x})= {\sum}_{\nu=1}^{N} {z_\nu}\,p({\|\bar{x}-\bar{x}_\nu\|}), \quad z_\nu \in
\mathbb{R}, \quad \bar{x}_\nu, \bar{x} \in \Omega=[-b,b]^3,
\end{equation}
is computationally demanding for large $N$. Since the generating radial basis function 
$p(\|\bar{x}\|)$ exhibits a slow polynomial decay in $1/\|\bar{x}\|$  as $\|\bar{x}\| \to \infty$, 
it follows that each individual term in (\ref{eqn:PotSum1}) contributes essentially to the total 
potential at each point in the computational domain $\Omega$. This predicts the $\mathcal{O}(N)$ 
complexity for a straightforward summation at every fixed space point $\bar{x} \in \mathbb{R}^3$. 
Moreover, in general, the radial function $p(\|\bar{x}\|)$ has a singularity or a cusp at the origin,
$\bar{x}=0$, making its accurate grid representation problematic. An efficient numerical scheme for 
the grid-based calculation of $P(\bar{x})$ in multiparticle systems can be constructed by using the 
RS tensor format \cite{BKK_RS:18}.

\subsection{Canonical tensor approximation of the Newton kernel}
 \label{ssec:Coulomb}

First, we recall the grid-based method for the low-rank canonical representation of a spherically 
symmetric kernel function $p(\|\bar{x}\|)$, $\bar{x}\in \mathbb{R}^d$ for $d=2,3,\ldots$, by its 
projection onto the set of piecewise constant basis functions, see \cite{BeHaKh:08} for the case of
the Newton kernel $p(\|\bar{x}\|)=\frac{1}{\|\bar{x}\|}$,  for $x\in \mathbb{R}^3$. A single 
reference potential like $1/\|\bar{x}\|$ can be represented on a fine 3D $n\times n\times n$ 
Cartesian grid as a low-rank canonical tensor \cite{HaKhtens:04I,BeHaKh:08}.

In the computational domain  $\Omega=[-b,b]^3$, let us introduce the uniform $n \times n \times n$ 
rectangular Cartesian grid $\Omega_{n}$ with mesh size $h=2b/n$ ($n$ even). Let $\{\psi_\textbf{i}\}$ 
be a set of tensor-product piecewise constant basis functions,
$ \psi_\textbf{i}(\bar{x})=\prod_{\ell=1}^3 \psi_{i_\ell}^{(\ell)}(\bar{x}_\ell)$,
for the $3$-tuple index ${\bf i}=(i_1,i_2,i_3)$, $i_\ell \in I_\ell=\{1,...,n\}$, $\ell=1,\, 2,\, 3 $.
The generating kernel $p(\|\bar{x}\|)$ is discretized by its projection onto the basis
set $\{ \psi_\textbf{i}\}$
in the form of a third order tensor of size $n\times n \times n$, defined entry-wise as
\begin{equation}  \label{galten}
\mathbf{P}:=[p_{\bf i}] \in \mathbb{R}^{n\times n \times n}, \quad
 p_{\bf i} =
\int_{\mathbb{R}^3} \psi_{\bf i} (\bar{x}) p(\|\bar{x}\|) \,\, \mathrm{d}{\bar{x}}.
\end{equation}
Then using the Laplace-Gauss transform and sinc-quadratures, the
$3$rd order tensor $\mathbf{P}$ can be approximated by
the $R$-term   canonical representation
(see \cite{HaKhtens:04I,BeHaKh:08,VeBoKh:Ewald:14} for details),
\begin{equation} \label{eqn:sinc_general}
    \mathbf{P} \approx  \mathbf{P}_R =
  \sum\limits_{k=1}^{R} {\bf p}^{(1)}_k \otimes {\bf p}^{(2)}_k \otimes {\bf p}^{(3)}_k
\in \mathbb{R}^{n\times n \times n}, 
\end{equation}
where ${\bf p}^{(\ell)}_k \in \mathbb{R}^n$.

Note that the reference tensor for summation of the potentials is generated in a
larger computational domain, $\widetilde{\bf P}_R  \in \mathbb{R}^{2n \times 2n \times 2n}$,
which is necessary for application of the shift-and-windowing transforms  ${\cal W}_\nu$,
see \Cref{ssec:RS_format} and  \cite{KhKhAn:12} for more details.

The canonical tensor representation of the Newton kernel (\ref{eqn:sinc_general}) has been 
successfully applied in computation of multidimensional operators in quantum chemistry
\cite{khor-ml-2009,KhKhAn:12}, where it was shown that calculations using the grid-based tensor 
approximations exhibit the same high accuracy level as the analytically based computation methods 
for the same multidimensional operators. The recent assembled tensor summation method of the 
long-range electrostatic potentials on large finite lattices \cite{VeBoKh:Ewald:14} has been proved 
to keep the rank of the collective potential on large 3D lattices to be as small as the rank of a 
canonical tensor for a single Newton kernel.

 \subsection{Properties of the RS tensor format} \label{ssec:RS_format}

Here we recall some properties of the  range separated (RS) tensor format introduced in 
\cite{BKK_RS:18} for modeling of the long-range interaction potential in multiparticle systems of 
general type. It is based on the partitioning of the reference tensor representation of the Newton 
kernel into long- and short-range parts with a following assembling of the collective electrostatic 
potential of a molecular system in a special way. According to the tensor canonical representation of 
the Newton kernel (\ref{eqn:sinc_general}) as a sum of Gaussians, one can distinguish their supports 
as  the short- and long-range parts,
 \[
  \mathbf{P}_R = \mathbf{P}_{R_s} + \mathbf{P}_{R_l},
\]
where
\begin{equation} \label{eqn:Split_Tens}
    \mathbf{P}_{R_s} =
\sum\limits_{k\in {\cal K}_s} {\bf p}^{(1)}_k \otimes {\bf p}^{(2)}_k \otimes {\bf p}^{(3)}_k,
\quad \mathbf{P}_{R_l} =
\sum\limits_{k\in {\cal K}_l} {\bf p}^{(1)}_k \otimes {\bf p}^{(2)}_k \otimes {\bf p}^{(3)}_k.
\end{equation}
Here, ${\cal K}_l := \{k|k = 0,1, \ldots, R_l\}$ and ${\cal K}_s := \{k|k = R_l+1, \ldots, M\}$ are 
the sets of indices for the long- and short-range canonical vectors. Then the optimal splitting 
(\ref{eqn:Split_Tens}) is applied to the reference canonical tensor ${\bf P}_R$ and to its 
accompanying version $\widetilde{\bf P}_R=[\widetilde{p}_R(i_1,i_2,i_3)]$,
$i_\ell \in \widetilde{I}_\ell$, $\ell=1,2,3$, such that
\[
 \widetilde{\bf P}_R = \widetilde{\mathbf{P}}_{R_s} + \widetilde{\mathbf{P}}_{R_l}
 \in \mathbb{R}^{2n \times 2n \times 2n}.
\]

The total electrostatic potential $P_0(\bar{x})$ in (\ref{eqn:PotSum1}) is represented by a projected 
tensor ${\bf P}_0\in \mathbb{R}^{n \times n \times n}$ that can
be constructed by a direct sum of shift-and-windowing transforms of the reference
tensor $\widetilde{\bf P}_R$ (see \cite{VeBoKh:Ewald:14} for more details),
\begin{equation}\label{eqn:Total_Sum}
 {\bf P}_0 = \sum_{\nu=1}^{N} {z_\nu}\, {\cal W}_\nu (\widetilde{\bf P}_R)=
 \sum_{\nu=1}^{N} {z_\nu} \, {\cal W}_\nu (\widetilde{\mathbf{P}}_{R_s} + \widetilde{\mathbf{P}}_{R_l})
 =: {\bf P}_s + {\bf P}_l.
\end{equation}

The shift-and-windowing transform ${\cal W}_\nu$ maps a reference tensor
$\widetilde{\bf P}_R\in \mathbb{R}^{2n \times 2n \times 2n}$ onto its sub-tensor
of smaller size $n \times n \times n$, obtained by first shifting the center of
the reference tensor $\widetilde{\bf P}_R$ to the grid-point $x_\nu$ and then restricting
(windowing) the result onto the computational grid $\Omega_n$.
However, the  tensor representation (\ref{eqn:Total_Sum}) is non-efficient as
the ranks are growing linearly in the number of particles and remain non-reducible 
in both canonical and Tucker tensor formats.

This problem is solved in \cite{BKK_RS:18} by considering  the global tensor
decomposition of only the  "long-range part" in the tensor ${\bf P}_0$, defined by
\begin{equation}\label{eqn:Long-Range_Sum} 
{\bf P}_l = \sum_{\nu=1}^{N} {z_\nu} \, {\cal W}_\nu (\widetilde{\mathbf{P}}_{R_l})=
\sum_{\nu=1}^{N} {z_\nu} \, {\cal W}_\nu
(\sum\limits_{k\in {\cal K}_l} \widetilde{\bf p}^{(1)}_k \otimes \widetilde{\bf p}^{(2)}_k
\otimes \widetilde{\bf p}^{(3)}_k).
\end{equation}

In \cite{BKK_RS:18} it was proven that the canonical rank of the tensor ${\bf P}_l$ only 
logarithmically depends on the number of particles $N$ involved in the summation. It was also shown 
that the rank reduction  ${\bf P}_l \mapsto {\bf P}_{R_L}$ may be efficiently implemented by using 
the C2T and T2C algorithms, where the RHOSVD is a key ingredient \cite{khor-ml-2009}. As for the 
short-range part of the collective potential, in the RS format it is represented by a single small 
size tensor supplemented by a list of the particles' coordinates.
  
According to the definition of the RS  canonical tensor format introduced in \cite{BKK_RS:18},
for the given separation parameter
$\gamma \in \mathbb{N}$, which is used for partitioning of the canonical vectors into
the long and short range contributions,
the RS-canonical tensor format specifies the class of $d$-tensors 
${\bf A}  \in \mathbb{R}^{n_1\times \cdots \times n_d}$
which can be represented as a sum of a   rank-${R_L}$ canonical tensor   
\begin{equation}
 {\bf A}_{R_L} = {\sum}_{k =1}^{R_L} \xi_k {\bf a}_k^{(1)} \otimes \cdots \otimes {\bf a}_k^{(d)}
\in \mathbb{R}^{n_1\times ... \times n_d}
 \end{equation}  
and
a (uniform) cumulated canonical tensor
\begin{equation}
  \widehat{\bf A}_S={\sum}_{\nu =1}^{N} c_\nu {\bf A}_\nu 
 \end{equation} 
  generated by shifts of a short reference tensor ${\bf A}_0$ with $\mbox{rank}({\bf A}_0)\leq R_0$.
 With these notations, the RS canonical tensor format reads as  
\begin{equation}\label{eqn:RS_Can}
 {\bf A} =  {\sum}_{k =1}^{R_L} \xi_k {\bf a}_k^{(1)}  \otimes \cdots \otimes {\bf a}_k^{(d)} +
 {\sum}_{\nu =1}^{N} c_\nu {\bf A}_\nu, 
\end{equation}
where   $\mbox{diam}(\mbox{supp}{\bf A}_\nu)\leq 2 \gamma$.

Lemma 3.9 in \cite{BKK_RS:18} presents the storage cost of the  RS-canonical
tensor ${\bf A}$ in (\ref{eqn:RS_Can}) as follows:
$$
\mbox{mem}({\bf A})\leq d R_L n + (d+1)N + d R_0 \gamma.
$$
Given ${\bf i}\in {\cal I} = I_1 \times \cdots \times I_d$, denote by
$\overline{\bf a}^{(\ell)}_{i_\ell}\in \mathbb{R_L}^{1\times R}$ the  
row-vector  with index $i_\ell$ in the side matrix $A^{(\ell)}\in \mathbb{R}^{n_\ell \times R_L}$
of ${\bf A}$, and let $\xi=(\xi_1,\ldots,\xi_d)$.
Then the ${\bf i}$-th entry of the RS-canonical tensor ${\bf A}=[a_{\bf i}]$ can be calculated
as a sum of long- and short-range contributions by
\[
 a_{\bf i}= \left(\odot_{\ell=1}^d \overline{\bf a}^{(\ell)}_{i_\ell} \right) \xi^T +
 \sum_{\nu \in {\cal L}({\bf i})} c_\nu {\bf A}_\nu({\bf i}),
\]
at the expense  $O(d R_L + 2 d \gamma R_0)$. 
Here, ${\cal L}({\bf i}) := \{\nu \in \{1,\ldots,N\}: {\bf i} \in \mbox{supp} {\bf A}_\nu\}$ is the set of 
indices which label all the short-range tensors ${\bf A}_\nu$ that include the grid point ${\bf i}$ 
within their effective support \cite{BKK_RS:18}. 
\begin{figure}[t]
\centering
\includegraphics[width=7.4cm]{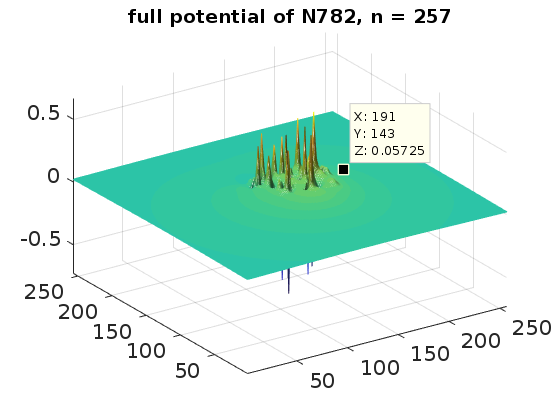}
\quad \quad
\includegraphics[width=7.4cm]{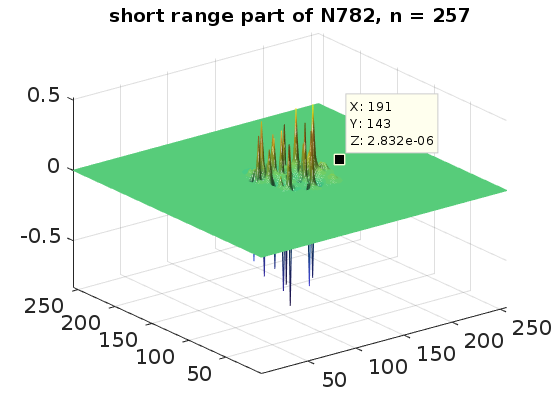}
\caption{The full free space electrostatic potential of a system with 782 particles (left)
and the sum of their short range contributions (right). }
\label{fig:Prot_1228}
\end{figure}

The RS tensor may be represented in a Tucker tensor format as well, see \cite{BKK_RS:18}.
RS tensors are efficient in many applications, for example for modeling of the electrostatics of
many-particle systems of general type, or for modeling of scattered multidimensional data
by using radial basis functions.
  
Next, we illustrate the performance of the canonical RS-tensor format in calculating the collective 
free space electrostatic potential of a model molecular system with 782 atoms, see also 
\cite{BKK_RS:18}. Summation is performed using the canonical tensor representing the reference
Newton potential computed on the $n\times n\times n$ 3D Cartesian grid with $n=257$ and the canonical 
rank $R=29$. The RS-tensor construction is performed with $r_\ell=14$, $R_s=15$, simply by dividing 
the canonical vectors into two groups, that is from every 29 vectors of the reference tensor and 
for every space dimension $(x,y,z)$, $15$ sharp Gaussians are separated as the short range 
part, and $14$ smoother Gaussians as the long-range 
part\footnote{Alternatively, separation of the reference tensor into the short- and long-range 
Gaussians may be performed by a chosen $\varepsilon$-truncation for the given interatomic distance.}. 
Then their contributions to the collective sum of potentials 
are calculated separately. To reduce the rank of the long-range collective sum,
the C2T (and T2C) transforms using RHOSVD \cite{khor-ml-2009} are used, with a choice of
the truncation threshold $10^{-8}$. The resulting rank 
decreases from $R_{\ell, total}=10948$, to $R_\ell=382$ 
and only logarithmically depends on the number of particles $N$ \cite{BKK_RS:18}.
\begin{figure}[t]
\centering
\includegraphics[width=7.0cm]{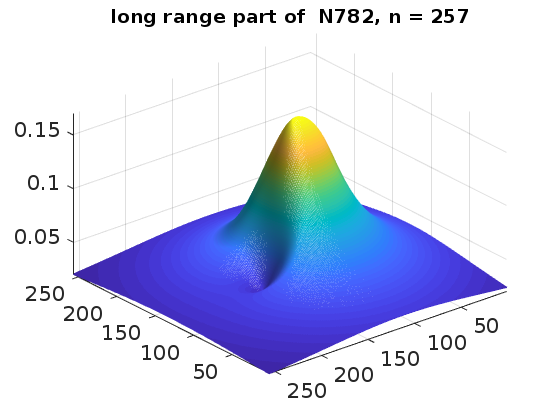}
\includegraphics[width=7.0cm]{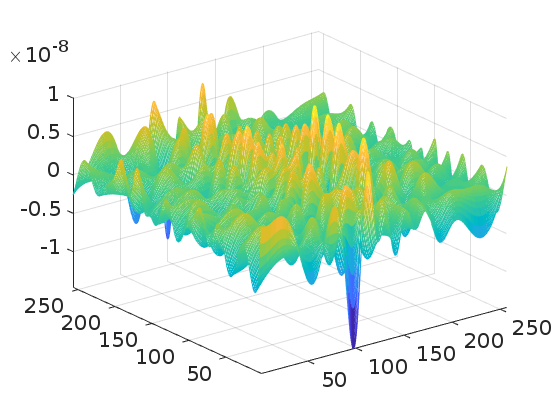}
\caption{The low-rank tensor representation of the long-range part in the
electrostatic potential of 782 charged particles (left)
and the error of the canonical rank reduction. }
\label{fig:Prot_long_1228}
\end{figure}

The left panel in Figure \ref{fig:Prot_1228} shows the cross-section of the collective electrostatic
potential of a molecular system at the middle of a $z$-plane, while the right panel shows
the cross-section of only the short-range part of the collective potential $\widehat{\bf A}_S$. 
Notice that in the left panel showing the total potential sum, the potential at the point of the 
plane with $(x,y)=(191,143)$ equals to $\sim 0.057$ units, while for the short range (right panel) 
the sum equals to $2.8\, 10^{-6}$. The left panel in Figure \ref{fig:Prot_long_1228} presents the 
cross-section of the low-rank long-range part of the collective potential at the same plane, while 
the right panel in this figure shows the error of the rank reduction (the same truncation threshold 
$10^{-8}$ as chosen above).

\section{Application of RS tensor format for solving PBE}\label{sec:RS_2_PBE}

The RS tensor formats can be gainfully applied in computational problems which include functions 
with multiple local singularities or cusps, Green kernels with essentially non-local behavior, as 
well as in various approximation problems  treated  by means of radial basis functions. In what 
follows, we describe the new approach  for the construction of computationally effective 
boundary/interface conditions and source terms in the Poisson-Boltzmann equation (PBE) describing 
the electrostatic potential of a biomolecule in gas phase and in solvent by solving the FEM/ FDM 
discretization of the regularized PBE. The main advantage of our approach is due to  complete 
avoidance of the direct FEM approximation (interpolation) of the highly singular right-hand sides in 
the traditional formulation of the PBE and, at the same time, preventing the modification of the 
stiffness matrix and/or the continuity conditions across the interface.

\subsection{Tensor based regularization scheme for the PBE }\label{ssec:PBE_Applic}

The traditional numerical approaches for solving the PBE are based on either multigrid  
\cite{Holst:2008} or domain decomposition \cite{CaMaSt:13} methods. Consider a solvated biomolecular 
system modeled by dielectrically separated domains with singular Coulomb potentials distributed in 
the molecular region. For schematic representation, we consider the system occupying a rectangular 
domain $\Omega$ with boundary $\partial \Omega$, see Figure \ref{fig:Protein}, where the solute 
(molecule) region is represented by $\Omega_m$ and the solvent region by $\Omega_s$, such that
\[
 \overline{\Omega} = \overline{\Omega}_m \cup \overline{\Omega}_s.
\]
The linearized Poisson-Boltzmann equation takes the form, see \cite{Holst:2008},
\begin{equation}\label{eqn:PBE}
 -\nabla\cdot(\epsilon \nabla u)+ \bar{\kappa}^2 u=\rho_f\quad \mbox{in } \quad \Omega,
\end{equation}
where $u$ denotes the target electrostatic potential of a protein, and 
$$
\rho_f= \sum\limits_{{k}=1}^N z_k \delta(\|\bar{x} - \bar{x}_k\|), \quad z_k \in \mathbb{R},
$$ 
is the scaled singular charge distribution supported at points $\bar{x}_k$ in $\Omega_m$, where
$\delta$ is the Dirac delta distribution, and $z_k\in \mathbb{R}$ denotes the charge located at 
the atomic center $\bar{x}_k$. Here $\epsilon_m=O(1)>0$ and $\bar{\kappa}=0$ in $\Omega_m$, while in 
the solvent region $\Omega_s$, we have $\bar{\kappa} \geq 0$ and $\epsilon_s \geq \epsilon_m$ 
(in some cases the ratio  $\epsilon_s/\epsilon_m$ could be about several tens).

The interface conditions on the interior boundary 
$\Gamma=\partial \Omega_m$ arise from the dielectric theory:
\begin{equation}\label{eqn:Intcond_PBE}
 [u]=0, \quad \left[ \epsilon\frac{\partial u}{\partial n}  \right]=0 \quad \mbox{on}\quad \Gamma.
\end{equation}
The boundary conditions on the external boundary
$\partial\Omega$ can be specified depending on the particular problem setting.
The simplest homogeneous Dirichlet boundary conditions $$u_{| \partial \Omega}=0$$ 
can be utilized. We shall also consider the practically interesting inhomogeneous Dirichlet boundary 
conditions taking the form
\begin{equation}\label{eqn:Dir_bc_PBE}
 u(\bar{x})_{| \partial \Omega}= \frac{1}{4 \pi \epsilon_m} \sum\limits_{{k}=1}^N 
 \frac{z_k e^{-\bar{\kappa} \|\bar{x} - \bar{x}_k\|}}{\|\bar{x} - \bar{x}_k\|}, 
 \quad \bar{x} \in \partial \Omega.
\end{equation}
The practically useful solution methods for the PBE are based on regularization 
schemes aiming at removing the singular component from the potentials in the governing equation.
Amongst others, we consider one of the most commonly used approaches based on the additive splitting 
of the potential in the molecular region $\Omega_m$, see for example \cite{Holst:2008}.
To that end, we first discuss   the additive splitting techniques introduced in  \cite{BKK_RS:16},
based on the application of the RS tensor format, 
\begin{equation}\label{eqn:regular_split}
 u=u^r +u^m_0, \quad \mbox{where}\quad u^m_0=0 \quad \mbox{in}\quad \Omega_s,
\end{equation}
and where the singular component $u^m_0$ satisfies the following Poisson equation in ${\Omega}_m$,
\begin{equation}\label{eqn:PE}
-\epsilon_m \Delta u^m_0 = \rho_f\quad \mbox{in}\quad  {\Omega}_m; 
\quad u^m_0=0 \quad \mbox{on}\quad {\Gamma}.
\end{equation}
In this case, equation (\ref{eqn:PBE}) can be transformed to an equation for the regular potential $u^r$:
\begin{equation}\label{eqn:Regul_PBE}
 -\nabla\cdot(\epsilon \nabla u^r)+ \bar{\kappa}^2 u^r=0 \quad \mbox{in } \quad \Omega,
\end{equation}
\[
 [u^r]=0,\quad 
 \left[ \epsilon\frac{\partial u^r}{\partial n}\right]=
 -\epsilon_m \frac{\partial u^m_0}{\partial n} 
 \quad \mbox{on}\quad \Gamma.
\]

To facilitate the solution of equation (\ref{eqn:PE}) with highly singular data in the right-hand side,
the singular potential $U$ in the free space was utilized, see \cite{BKK_RS:18},
\begin{equation}\label{eqn:Free_spa_Pot}
-\epsilon_m \Delta U = \rho_f \quad  \mbox{in}\quad \mathbb{R}^3,\quad |U(\bar{x})|\to 0, 
\quad |\bar{x}|\to \infty,
\end{equation}
that can be written in the explicit form
\[
 U(\bar{x})= \frac{1}{4\pi \epsilon_m} \sum\limits_{{k}=1}^N \frac{z_k}{\|\bar{x} - \bar{x}_k\|}.
 \]
Introduce the characteristic (indicator) function, $\chi_{[\Omega_m]}(\bar{x})$, $\bar{x} \in \Omega$, 
of the domain $\Omega_m\subset \Omega$ by
\begin{equation}\label{eqn:Charact_func}
 \chi_{[\Omega_m]}(\bar{x}) =
  \begin{cases}
   1 & \mbox{if } \bar{x} \in \overline{\Omega}_m\\
   0 & \mbox{if } \bar{x} \in \Omega_s=\Omega \setminus \overline{\Omega}_m.
  \end{cases}
\end{equation}
Then the restriction of $U$ onto $\Omega_m$ can be calculated by 
$$u^m=\chi_{[\Omega_m]} U, 
$$
implying the decomposition 
$$
u^m_0=u^m + u^{harm},
$$
where the harmonic function $u^{harm}$ compensates the discontinuity of $u^m$ on $\Gamma$,
\begin{equation*}\label{eqn:Harmonic}
 \Delta u^{harm} = 0 \quad \mbox{in}\quad  {\Omega}_m; \quad u^{harm}=-u^m=-U \quad\mbox{on}\quad \Gamma.
 \end{equation*}

The advantage of this formulation is twofold: 
\begin{enumerate}[label=(\alph*)]
 \item the absence of singularities in the solution $u^r$, and
 \item the localization of the solution splitting only in the domain $\Omega_m$.
\end{enumerate}
Grid representation of the free space singular potential $U$, which may include a sum of hundreds 
or even thousands of single Newton kernels in 3D, leads to a challenging computational problem.
In our approach it can be represented on large tensor grids in $\Omega$ with controlled 
precision by using the RS tensor format  \cite{BKK_RS:18}
characterized by the separability constant $\gamma >0$ which in our application 
can be associated with the
Van Der Waals inter-atomic distance, see \Cref{ssec:RS_format}. 
The long-range component in the formatted parametrization
remains smooth and allows global low-rank representation in $\Omega$. 
We conclude with the following 
\begin{proposition}\label{pro:PBE_RSTensor}
Let the effective support of the short-range components in the reference
 potential ${\bf P}_R$ be chosen  not larger than $\gamma/2$. Then 
 the interface conditions in the regularized formulation of the PBE in (\ref{eqn:Regul_PBE})
 depend only on the low-rank long-range component in the free-space electrostatic potential
 of the atomic system. The numerical cost to  build up the interface conditions 
 on $\Gamma$ in (\ref{eqn:Regul_PBE})  does not depend on the number of particles $N$.
\end{proposition}

Notice that here, we describe the splitting scheme in (\ref{eqn:PE}) -- (\ref{eqn:Regul_PBE})
just for illustration of the applicability of the RS 
tensor format for the solution of the PBE. This scheme requires modification of the interface conditions that 
is equivalent to a change of the system matrix which leads to a complicated implementation scheme. To 
avoid this nontrivial task, in what follows, we introduce an alternative approach, which will be 
discussed in the next section.

The regularization $u=u^r+u^m_0$  like in (\ref{eqn:PE}) -- (\ref{eqn:Regul_PBE})  
benefits from the local-global separability in the low-rank RS tensor representation of the 
free space electrostatic potential. 
Notice that here, we describe the splitting scheme in (\ref{eqn:PE}) -- (\ref{eqn:Regul_PBE})
just for illustration of the applicability of the RS 
tensor format for the solution of the PBE. This scheme requires modification of the interface conditions that 
is equivalent to a change of the system matrix which leads to a complicated implementation scheme. To 
avoid this nontrivial task, in what follows, we introduce an alternative approach, 
which avoids the additional computation 
of the auxiliary harmonic function $u^{harm}$ in the rather complicated domain $\Omega_m$ as well as 
the modification of the interface conditions.

\subsection{The new RS tensor based splitting scheme} \label{ssec:Disc_Split_scheme}

In this section, we present the new splitting scheme which is based on the range separated 
representation of the Dirac $\delta$-distribution \cite{khor-DiracRS:2018}, 
which composes the highly singular right-hand side in the target PBE (\ref{eqn:PBE}) 
or   Poisson  equation (PE) (\ref{eqn:PE}). Here, we consider the PE as proof of concept and validate 
the numerical results in \Cref{sec:Numer_Tests}. The idea is to modify the right-hand side $\rho_f$
in such a way that the short-range part in the solution $u$ can be pre-computed independently 
by the direct tensor decomposition of the free space potential, and the initial elliptic equation applies 
only to the long-range part of the total potential, satisfying the equation with the modified 
right-hand side by using the RS splitting of the Dirac delta, see \cite{khor-DiracRS:2018} for more details. 
The latter is a smooth function, hence the FDM/FEM 
approximation error can be reduced dramatically even on relatively coarse grids in 3D.

For ease of presentation, we consider the simplest case of the single atom with unit charge located at the 
origin, such that the exact electrostatic potential reads $u(\bar{x})=\frac{1}{\|x\|}$, $x\in \mathbb{R}^3$. 
Recall that the Newton kernel (\ref{eqn:sinc_general}) discretized by the $R$-term sum 
of Gaussian type functions 
living on the $n \times n \times n$ tensor grid $\Omega_n$ is represented by a sum of 
short- and long-range tensors,
\[
 \frac{1}{\|x\|}   \rightsquigarrow \mathbf{P}_R = \mathbf{P}_{R_s} + \mathbf{P}_{R_l} 
 \in \mathbb{R}^{n \times n \times n},
\]
where $\mathbf{P}_{R_s}$ and $\mathbf{P}_{R_l}$ are defined in (\ref{eqn:Split_Tens}).

Let us formally discretize  the exact equation for the Newton potential, $u(\bar{x})=\frac{1}{\|x\|}$,
\[
 -\Delta \frac{1}{\|x\|} =4 \pi \,\delta(\bar{x}),
\]
by using the FDM/FEM Laplacian matrix $A_{\Delta}$ instead of $\Delta$ and via substitution of 
the canonical tensor decomposition $\mathbf{P}_R $ instead of $u(\bar{x})=\frac{1}{\|x\|}$.
This leads to the grid representation  of the discretized Dirac delta \cite{khor-DiracRS:2018}
\[
 \delta(\bar{x}) \rightsquigarrow \boldsymbol{\delta}_h:=-\frac{1}{4 \pi}A_{\Delta} \mathbf{P}_R,
\]
that will be applied in the framework of our discretization scheme. 

We remind that the 3D finite difference Laplacian matrix $A_{\Delta}$, 
defined on the uniform rectangular grid, takes the form
\begin{equation}\label{eqn:Lapl_Kron3}
A_{\Delta} = \Delta_{1} \otimes I_2\otimes I_3 + I_1 \otimes \Delta_{2} \otimes I_3 + 
I_1 \otimes I_2\otimes \Delta_{3},
\end{equation}
where 
$-\Delta_\ell = h_\ell^{-2} \mathrm{tridiag} \{ 1,-2,1 \} 
\in \mathbb{R}^{n_\ell \times n_\ell}$, $\ell=1,2,3$, 
denotes the discrete univariate Laplacian, such that the Kronecker rank 
 of $A_{\Delta}$ equals to $3$. Here $I_\ell$, $\ell=1,2,3$,  
 is the identity matrix in the corresponding single dimension. 
 
Now we are in the position to describe the RS tensor based splitting scheme.
To that end, we use the splitting of the discretized $\delta$-distribution 
into short- and long-range components in the form \cite{khor-DiracRS:2018},
\begin{equation}\label{Delta_sum}
 \boldsymbol{\delta}_h= \boldsymbol{\delta}_s + \boldsymbol{\delta}_l,
\end{equation}
where
\begin{equation}\label{LS_delta}
 \boldsymbol{\delta}_s:= -\frac{1}{4 \pi}A_{\Delta} \mathbf{P}_{R_s}, \quad 
 \mbox{and} \quad \boldsymbol{\delta}_l:= -\frac{1}{4 \pi}A_{\Delta} \mathbf{P}_{R_l}.
\end{equation}
We recall that by construction, the short range part vanishes on the interface $\Gamma$,
hence it satisfies the discrete Poisson equation in $\Omega_m$ with the respective right-hand
side in the form $\boldsymbol{\delta}_s$ and zero boundary conditions on $\Gamma$ \cite{khor-DiracRS:2018}. 
Then we deduce that this equation can be subtracted from the full
discrete linear system, such that the long-range component of the solution, $\mathbf{P}_{R_l}$,
will satisfy the same linear system of equations (same interface conditions), 
but with a modified right-hand side 
corresponding to the weighted sum of the long-range tensors  $\boldsymbol{\delta}_l$ only.
In the simple example of the single charge, we arrive at the particular 
discrete Poisson equation for the long-range part in the full potential 
$\mathbf{P}_R$, ${\bf U}_l=\mathbf{P}_{R_l}$, 
\begin{equation}\label{delta_L}
 - A_{\Delta} {\bf U}_l = \boldsymbol{\delta}_l,
\end{equation}
which can be solved by an appropriate iterative method.

\begin{figure}[t]
\centering
\includegraphics[width=8.1cm]{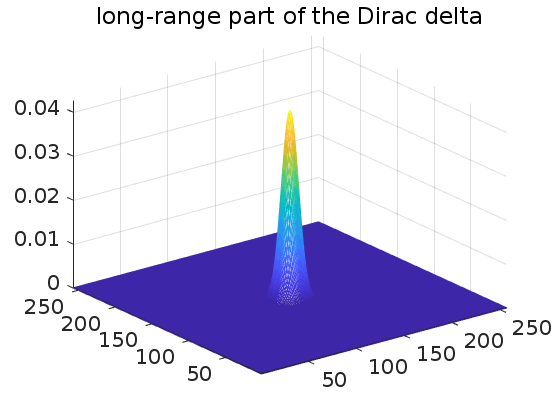} 
\caption{The long-range part of the Dirac delta $\boldsymbol{\delta}_l$ on an 
$n^{\otimes 3}$ 3D grid, $n=256$. }
\label{fig:Single_RHS_delta}
\end{figure}

Figure \ref{fig:Single_RHS_delta} illustrates the modified 
right hand side representing the long-range part of the discrete 
Dirac delta $\boldsymbol{\delta}_l$.
It it worth noting that the FEM approximation theory can be applied to this formulation since 
the input data (i.e., the right-hand side) are regular enough. 
However, it is not the case for the initial formulation with 
the highly singular Dirac delta distribution in the right-hand side. 

This scheme can be easily extended to the case of many-atomic systems just by additive representation of 
the short- and long-range parts in the total free space potential, 
\begin{equation}\label{delta_RL}
 - A_{\Delta} {\bf P}_{R_l} = \boldsymbol{\delta}_{R_L},
\end{equation}
where we suppose that $R_L$ is the rank of the long-range part ${\bf P}_{R_l}$
of the corresponding RS tensor of type (\ref{eqn:RS_Can}), and 
$\boldsymbol{\delta}_{R_L}$ is calculated as shown in (\ref{tensor_delta}). 

We summarize the following benefits of the aforementioned solution decomposition scheme. 
 \begin{itemize}
 \item Most important is that due to  efficient  splitting of the short- and 
long-range parts in the target tensor representing both the single Newton kernel 
and the total free-space potential, there is no need to modify jump conditions at the interface.
\item A remarkable advantage is that the long-range part in the RS tensor decomposition 
of the Dirac-delta distribution \cite{khor-DiracRS:2018}
vanishes at the interface and, hence, the modified right-hand side generated by this 
long-range component remains localized in the ``linear" solute region.
\item The boundary conditions are obtained from the long-range part in the tensor 
representation of the collective electrostatic potential which reduces the computational 
costs involved in solving some external analytical function at the boundary.
\item Only a single system of algebraic equations is solved for the smooth long-range 
(i.e., regularized) part of the collective potential discretized with controllable precision
on a relatively coarse grid, which is then added to the directly 
precomputed (avoiding PDE solutions) low-rank tensor representation for the short-range 
contribution.
\end{itemize}

Next, we briefly comment on the general FEM approximation for the Laplacian.
We recall the tensor-based scheme for evaluation of the Laplace operator 
in a separable basis set  \cite{KhKhAn:12} 
applied for calculation of the kinetic energy part in the Fock operator. 
  
Let the problem be posed in the finite volume box $ \Omega=[-b,b]^3 \in \mathbb{R}^3$,
subject to the homogeneous Dirichlet boundary conditions on $\partial \Omega$.
For given discretization parameter $n \in \mathbb{N}$,   
the equidistant $n\times n \times n$ tensor grid $\omega_{{\bf 3},N}=\{x_{\bf i}\} $, 
${\bf i} \in {\cal I} :=\{1,...,n\}^3 $ is used, with the mesh-size $h=2b/(n +1)$. Then the Kronecker 
rank-$3$ tensor representation of the respective Galerkin FEM stiffness matrix is 
given by (\ref{eqn:Lapl_Kron3}).  

Notice that the MATLAB representation of the matrix $A_\Delta$ 
(say, the FD matrix) can be easily described in terms of $\texttt{kron}$ operations as follows
\begin{equation}\label{3D_Lap}
 \frac{1}{h^2}A_\Delta = \texttt{kron}(\texttt{kron}(\Delta_1,I ),I ) + 
 \texttt{kron}(\texttt{kron}(I ,\Delta_1,I )+\texttt{kron}(\texttt{kron}(I ,I ),\Delta_1),
\end{equation}
applied to a long vector of size $n^3$ representing the Newton potential.

Then the rank-structured calculation of the ``collective'' right-hand 
side $\boldsymbol{\delta}_{R_L}$ 
 in (\ref{delta_RL}) is reduced to one-dimensional operations, 
\begin{equation}\label{tensor_delta}
 -\boldsymbol{\delta}_{R_L} = 
 \sum\limits_{k=1}^{R_L} \xi_k ( \Delta_1 {\bf a}^{(1)}_k \otimes {\bf a}^{(2)}_k \otimes {\bf a}^{(3)}_k 
 + {\bf a}^{(1)}_k \otimes \Delta_1 {\bf a}^{(2)}_k \otimes {\bf a}^{(3)}_k 
 +  {\bf a}^{(1)}_k \otimes  {\bf a}^{(2)}_k \otimes \Delta_1 {\bf a}^{(3)}_k), 
\end{equation}
where ${\bf a}^{(\ell)}_k$, $\ell=1,2,3$, are the canonical vectors and $R_L$ is the canonical rank 
of the long-range part of the collective electrostatic free space potential of a biomolecule 
computed in the RS tensor format (\ref{eqn:RS_Can}). 

This is the tensor ansatz to be used as the right-hand side in the equation (\ref{delta_L}),
which we apply in numerical experiments.
With a subsequent usage of the canonical-to-full tensor transform and after reshaping a three-tensor
into a long vector, $\boldsymbol{\delta}_{R_L}$ 
is applied in a standard PBE iterative solver as the RHS for the long-range part.
Another advantage of our scheme is that the short-range part of the solution in the 
PBE (\ref{eq:PBE}) is obtained for free, since it is 
merely incorporated as the set of short-range parts of the 
respective Newton potentials for every particle in a biomolecule. That corresponds to a set of
tensors in the second term of the collective electrostatic potential in the
RS tensor format (\ref{eqn:RS_Can}).

\subsection{Discussion of the computational scheme}\label{ssec:Comp_scheme}

We summarize the main computational tasks involved in the presented tensor-based numerical scheme. The 
regularized PBE (RPBE) can be constructed as follows:
\begin{enumerate}[label=(\alph*)]
 \item Compute the regularized component $\boldsymbol{\delta}_{R_L}$ of the discretized Dirac delta  
 distribution in (\ref{delta_RL}) as described in \cite{khor-DiracRS:2018};
 \item Substitute $\boldsymbol{\delta}_{R_L}$ into the right-hand side of the PBE (\ref{eqn:PBE}) to 
 obtain the following RPBE
 \begin{equation}\label{eqn:RPBE_decomp}
-\nabla \cdot(\epsilon\nabla u^r(\bar{x})) + \bar{\kappa_2}^2(\bar{x})u^r(\bar{x}) 
   = \boldsymbol{\delta}_{R_L},   \quad \mbox{in} \,\, \Omega,
\end{equation}
subject to the boundary condition  in (\ref{eq:DH_solution}).
 \item Discretize the RPBE (\ref{eqn:RPBE_decomp}) to obtain the following linear system of equations
\begin{equation}\label{eqn:RPBE_system}
Au^r = b,
\end{equation}
which can be solved by any linear system solver.
\item Obtain the final PBE solution $u$ by the sum 
\[u = u^r + u^s,\] where $u^s = \mathbf{P}_{R_s}$ is the precomputed short-range component of the 
Newton potential sum, see (\ref{eqn:Split_Tens}).
\end{enumerate}

The splitting scheme described in \Cref{ssec:Disc_Split_scheme} allows to reduce 
the initial equation to the solution of the system with 
modified right-hand side by using the  range-separated representation of the discretized 
Dirac delta. The problem is reduced to 
the direct tensor-based computation (without solving a PDE) of the short-range 
part in the collective free space electrostatic potential, see  (\ref{eqn:Split_Tens}), 
and to the 
subsequent solution of the PBE equation for the long-range part only by the simple modification of the
right-hand side. The advantage is that the PBE applies to the smooth part, $u^r$, of the total potential 
and hence a controllable FDM/FEM approximation error on moderate size 3D grids can be achieved.

This method  can be combined with the reduced basis approach for PBE with parametric 
coefficients to further accelerate the numerical computations \cite{KwHeFeStBe:16,KwBFM2017}. This is 
because the modified model is affinely dependent on the parameter $\bar{\kappa}$, 
thereby providing a natural off-line/on-line decomposition of the reduced basis method.

Finally, we notice that the important characterization of the protein molecule is
given by the electrostatic solvation energy \cite{Holst:2008}, which is the difference 
between the electrostatic free energy in the solvated state (described by the PBE) 
and the electrostatic free energy in the absence of solvent. 
Having at hand the free energy $E_N$ (see \cite{BKK_RS:16}, Lemma 4.1),  the electrostatic solvation energy 
can be computed in the framework of the regularized formulation of the PBE as described above.

\section{Numerical Tests}
\label{sec:Numer_Tests}

In this section, we consider only the free space electrostatic potential for the modified PE and the 
RS tensor format based splitting scheme. We compare the results with those of the 
traditional PE for various biomolecules. Notice that the PBE can be reduced to the PE by considering 
the case of zero ionic strength which implies that the function $\bar{\kappa}^2(\bar{x})$, hence the Boltzmann 
distribution term, in (\ref{eq:PBE}) is annihilated. Consequently, homogeneous dielectric constants of 
$\epsilon_m = \epsilon_s = 1$ are considered. We compute the electrostatic potentials using 
$n\times n\times n$ 3D Cartesian grids, in a box $[-b,b]^3$ with equal step size $ h =\frac{2b }{n-1}$. 
Note that computations by the PBE/PE solver are bounded by $n=257$, on a PC with $8$GB RAM due to the 
storage needs of the order of $\mathcal{O}(n^3)$. Note that since the storage requirements for the 
tensor-structured representation of the electrostatic potential is of the order of $\mathcal{O}(n)$, 
see (\ref{eqn:sinc_general}), in tensor computations of molecular potentials one can use 
much finer 3D grids.

\begin{figure}[t]
\centering
\includegraphics[width=5.1cm]{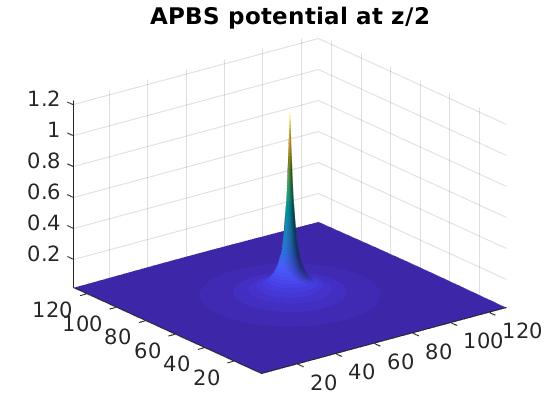} 
\includegraphics[width=5.0cm]{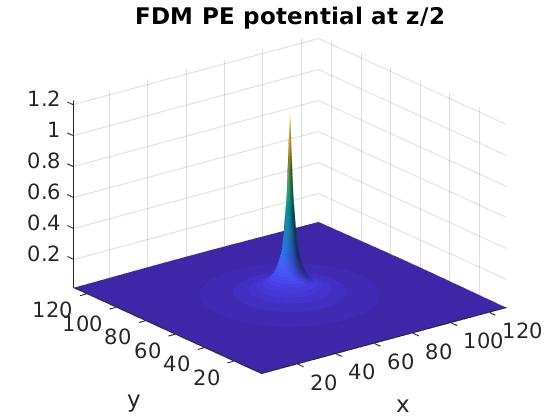} 
\includegraphics[width=5.3cm]{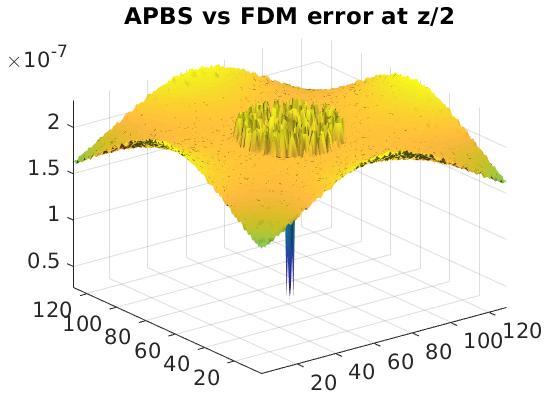}
\caption{The free space potential for the Born ion computed by the APBS (left), the FDM solver (middle) 
and the corresponding error (right). }
\label{fig:APBS_FDM_error}
\end{figure}

First, we validate the FDM solver for the PE by comparing its solution with that of the adaptive 
Poisson-Boltzmann solver (APBS) using the finite difference multigrid calculations with PMG \footnote{PMG 
is a Parallel algebraic Multigrid code for General semilinear 
elliptic equations.} option \cite{Bakersept2001, Chen:07}. Here, we consider a Born ion of 
unit charge, unit ($\textrm{\AA}$) radius and located at the origin ($(0,0,0)$). Figure 
\ref{fig:APBS_FDM_error} shows the single Newton kernel for the Born ion, approximated by the PE on 
an $n\times n$ grid surface with $n=129$ at the cross section of the volume box in the middle of the 
z-axis computed by the APBS and the FDM solvers and the corresponding error between the two solutions. The 
results show that the FDM solver provides as accurate results as those of the APBS with a discrete $L_2$ 
error of $\mathcal{O}(10^{-9})$ in the full solution. Similar results for 
varying proteins are illustrated in \cite{KwBFM2017}. 

\begin{figure}[b]
\centering
\includegraphics[width=5.1cm]{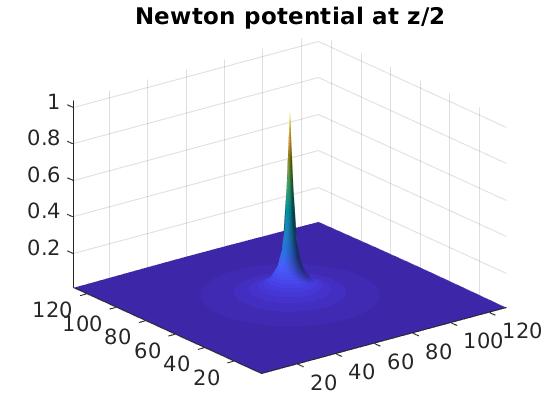} 
\includegraphics[width=5.0cm]{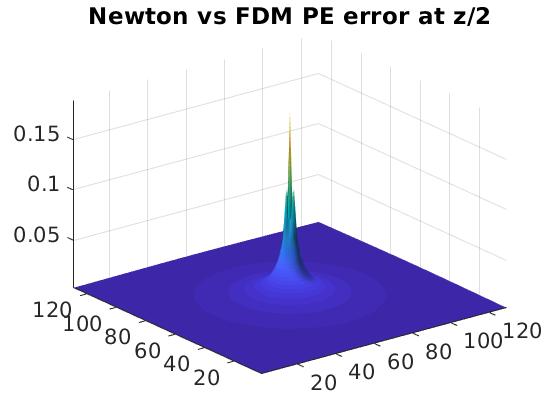} 
\includegraphics[width=5.3cm]{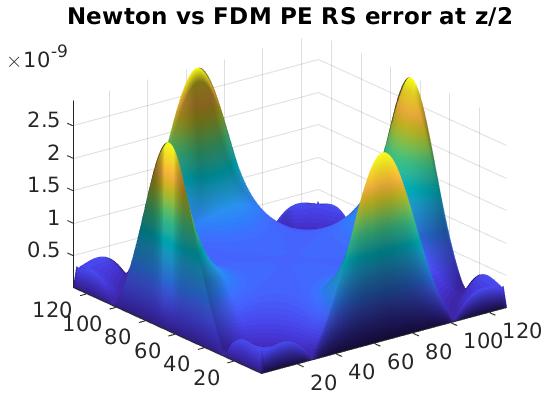}
\caption{The Newton potential computed by the canonical tensor decomposition (left), the error of its 
computation on the same grid by using the classical PE (middle) and by the modified PE (right). }
\label{fig:Newton_error}
\end{figure}

Secondly, we compare the accuracy of the traditional PE model and of the PE model modified by the RS 
tensor format for the approximation of the single Newton kernel. Figure \ref{fig:Newton_error} shows the 
single Newton kernel on an $n\times n$ grid surface with $n=129$ at the cross section of the volume box in 
the middle of the z-axis computed by the canonical tensor approximation obtained by sinc-quadratures and 
the corresponding errors by the traditional PE model and the modified PE model computed by the FDM solver. 
We notice that the solution of the modified PE model is of higher accuracy than that of the traditional 
PE model because it captures the singularities exactly, see the central region of \Cref{fig:Newton_error} 
(right).

\begin{figure}[t]
\centering
\includegraphics[width=5.1cm]{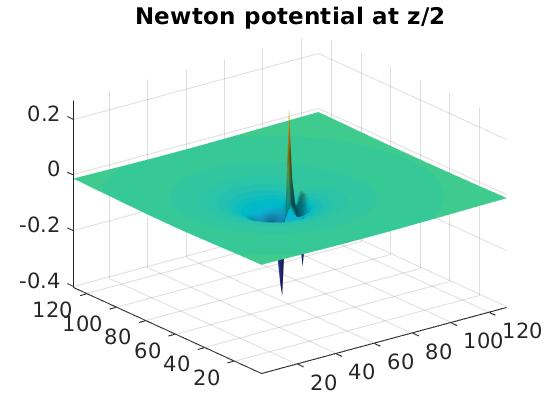} 
\includegraphics[width=5.0cm]{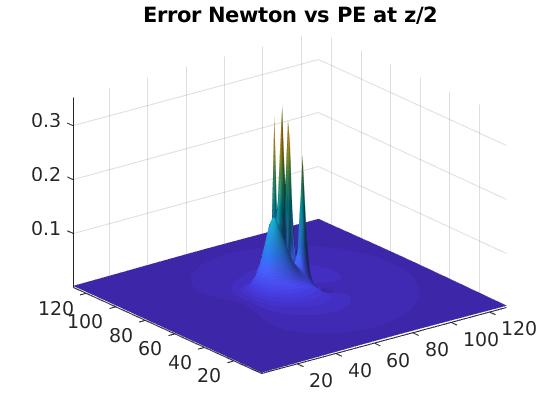} 
\includegraphics[width=5.3cm]{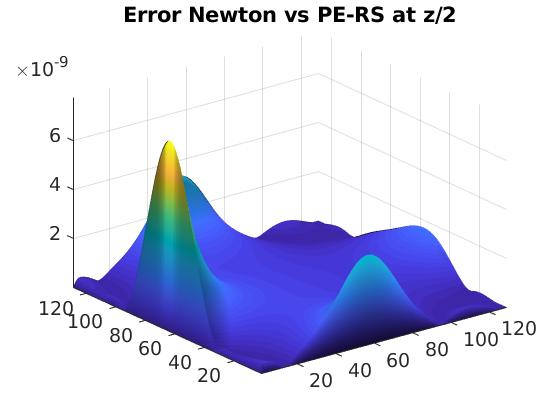}
\caption{\small The Newton potential sums computed by the canonical tensor decomposition (left), the error 
of its computation on the same grid by using the classical PE (middle) and by the modified PE (right).}
\label{fig:Pot_18_error}
\end{figure}

In a similar vein, we consider the acetazolamide compound consisting of $18$ atoms and determine the 
accuracy of the traditional PE model vis a vis the PE model modified by the RS tensor format.
Figure \ref{fig:Pot_18_error} shows the sum of the electrostatic potentials for the acetazolamide, 
computed on an $n\times n$ grid surface with $n=129$ by the canonical tensor approximation obtained by 
sinc-quadratures and the corresponding errors by the traditional PE model and the modified PE model 
computed by the FDM solver. It is clearly shown that the modified PE model provides highly accurate 
solutions as compared to those of the traditional PE due to the accurate treatment of the solution 
singularities by the RS tensor format.   

\begin{figure}[t]
\centering
\includegraphics[width=5.1cm]{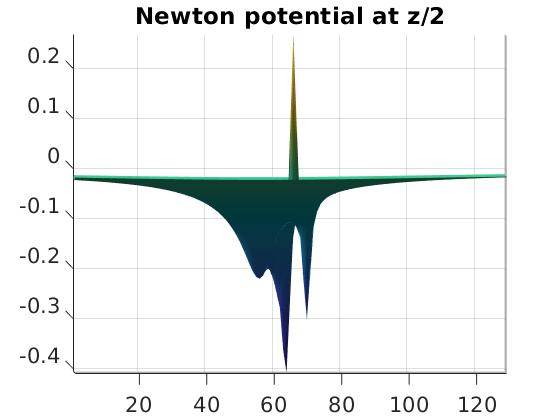} 
\includegraphics[width=5.1cm]{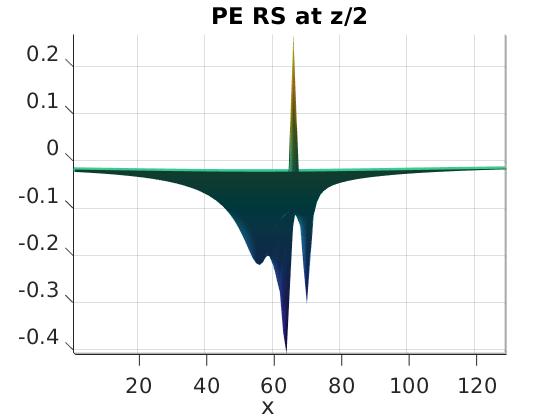}
\includegraphics[width=5.1cm]{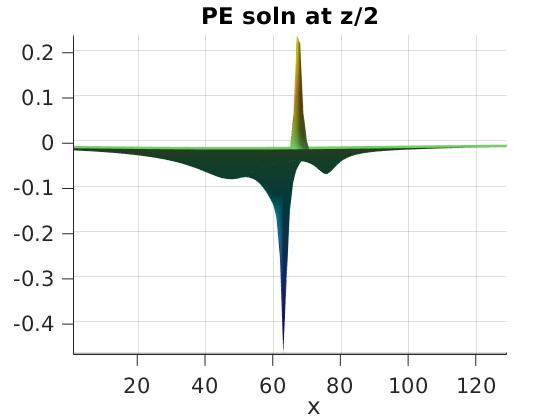}
\caption{\small Demonstration of the solution singularities for the acetazolamide molecule captured 
by the canonical tensor approximation (left), by the modified PE model (middle) and by the classical 
PE (right).}\label{fig:Pot_18_sing}
\end{figure}

We notice, in addition, that the classical PE model does not capture accurately the singularities in 
the electrostatic potential due to the numerical errors introduced by the Dirac delta 
distribution and partly due to the smoothing effect caused by the spline interpolation of the charges 
onto the grid. The modified PE model, on the other hand, is able to capture the singularities due to the 
independent treatment of the singularities by the RS tensor technique. This is clearly demonstrated in 
Figure \ref{fig:Pot_18_sing}.

\subsection{Solutions to the modified Poisson equation on a sequence of fine grids}
\label{ssec:Poisson_RS} 
 
\begin{figure}[t]
\centering
\resizebox{17cm}{7.0cm}{%
\begin{tikzpicture}

\node (Newt_full_257) at (0,3.0) {\begin{tabular}{l}
  {\includegraphics[width=0.48\textwidth]{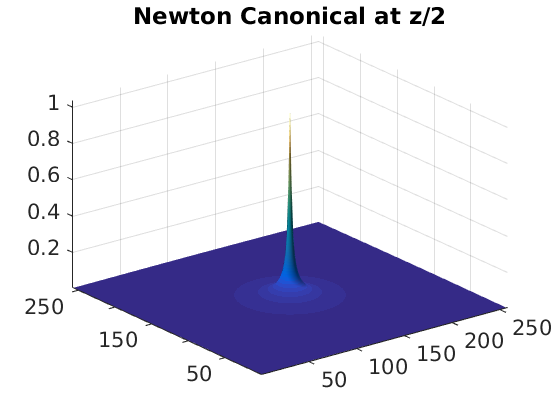}}
 \end{tabular}
 };
 
 \node (PBE_RS_full_257) at (0,-3.0) {\begin{tabular}{l}
  {\includegraphics[width=0.48\textwidth]{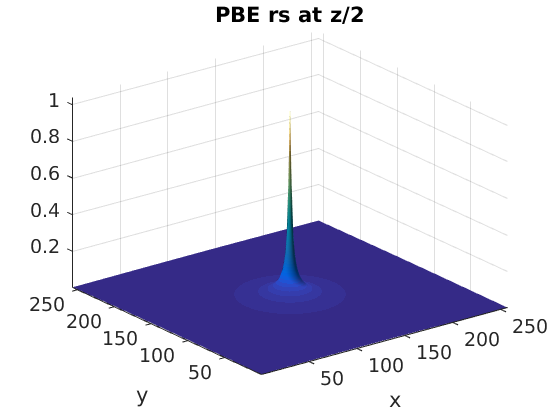}}
 \end{tabular}
 };
 
 \node (Error_PBE_RS_Newt_257) at (10,0) {\begin{tabular}{l}
  {\includegraphics[width=0.48\textwidth]{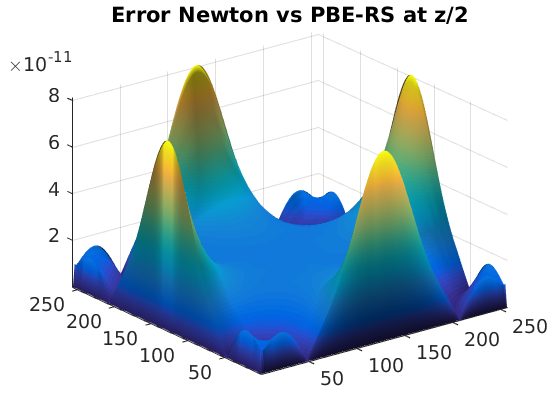}}
 \end{tabular}
 };
 \draw[->,draw=blue,thick] (3.5,0) -- (Error_PBE_RS_Newt_257.west);
 
 \end{tikzpicture}}
 \caption{\label{fig:Newt_pot_err_257}
 Absolute error between the solutions of the Newton potential and the modified Poisson equation for the 
 Born ion.}
\end{figure}

Here, we illustrate the accuracy of the modified PE by calculating the free space electrostatic potential 
on a sequence of fine grids and compare with the solution of the exact Newton potential determined by the 
canonical tensor representation. We first consider the aforementioned Born ion case and show the absolute 
error for the finest Cartesian grid and the discrete $L_2$ norm of the error for a sequence of Cartesian 
grids. Figure \ref{fig:Newt_pot_err_257} shows the absolute error of $\mathcal{O}(10^{-11})$ obtained on 
a $257^3$ Cartesian grid and $32\,\textrm{\AA}$ box length. Table \ref{Tab:Elect_error_Born} shows the 
decay of the discrete $L_2$ norm of the error for a sequence of grid refinements.
\begin{table}[t]
  \begin{center}%
  \begin{tabular}
  [c]{|c|c|c|c|}%
  \hline
    n         & 97   & 129   & 257   \\
    \hline
  Discrete $L_2$ norm  & $4.1176\times10^{-7}$   & $9.5516\times10^{-8}$ & $2.7975\times10^{-9}$  \\
  \hline
  \end{tabular}
  \caption{\small The discrete $L_2$ norm of the error with respect to grid size for the Born ion.}
  \label{Tab:Elect_error_Born}
  \end{center}
\end{table}

\begin{figure}[H]
\centering
\resizebox{17cm}{7.0cm}{%
\begin{tikzpicture}

\node (Acet_Newt_full_257) at (0,3.0) {\begin{tabular}{l}
  {\includegraphics[width=0.48\textwidth]{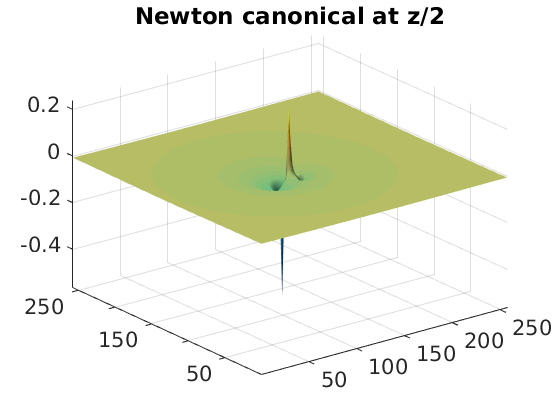}}
 \end{tabular}
 };
 
 \node (PBE_RS_full_Acet_257) at (0,-3.0) {\begin{tabular}{l}
  {\includegraphics[width=0.48\textwidth]{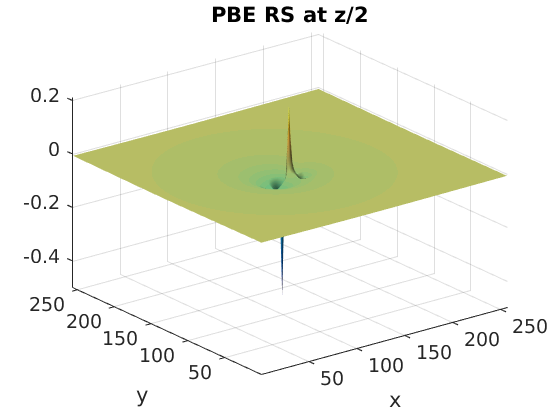}}
 \end{tabular}
 };
 
 \node (Error_PBE_RS_Newt_Acet_257) at (10,0) {\begin{tabular}{l}
  {\includegraphics[width=0.48\textwidth]{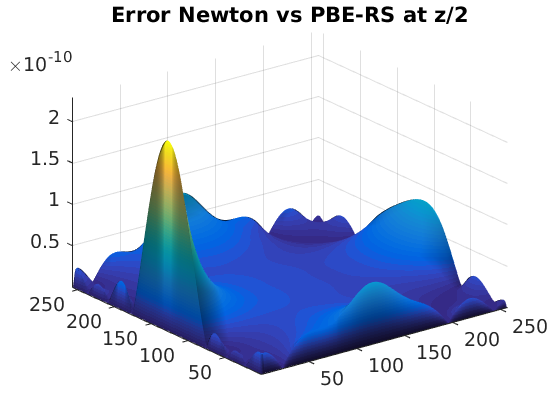}}
 \end{tabular}
 };
 \draw[->,draw=blue,thick] (3.5,0) -- (Error_PBE_RS_Newt_Acet_257.west);
 
 \end{tikzpicture}}
 \caption{\label{fig:Acet_Newt_pot_err_257}
 Absolute error between the solutions of the Newton potential sums and the modified Poisson equation
 for the acetazolamide molecule.}
\end{figure}

Next, we consider the acetazolamide molecule with $18$ atoms. This molecule is used as a ligand 
in the human carbonic anhydrase (hca) protein-ligand complex for the calculation of the binding 
energy in the adaptive Poisson-Boltzmann software (APBS) package \cite{Holst:93} and the MATLAB 
program for biomolecular electrostatic calculations, (MPBEC) \cite{Vergara-Perez2016}.  
The electrostatic  potential is computed as in the previous case, employing the same grid properties. 
The absolute error is shown in Figure \ref{fig:Acet_Newt_pot_err_257} and the error behaviour with 
respect to mesh refinements is shown in Table \ref{Tab:Elect_error_Acet}.

\begin{figure}[t]
\centering
\resizebox{17cm}{7.0cm}{%
\begin{tikzpicture}

\node (379_Newt_full_257) at (0,3.0) {\begin{tabular}{l}
  {\includegraphics[width=0.48\textwidth]{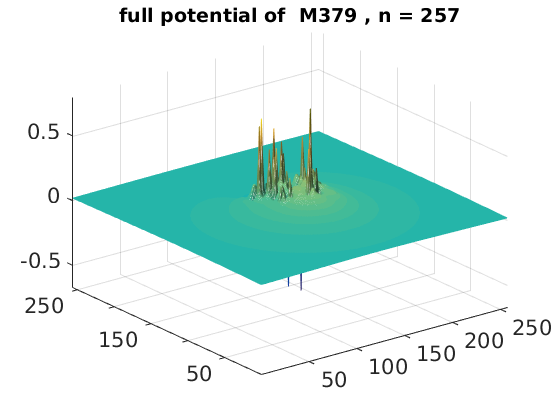}}
 \end{tabular}
 };
 
 \node (PBE_RS_full_379_257) at (0,-3.0) {\begin{tabular}{l}
  {\includegraphics[width=0.48\textwidth]{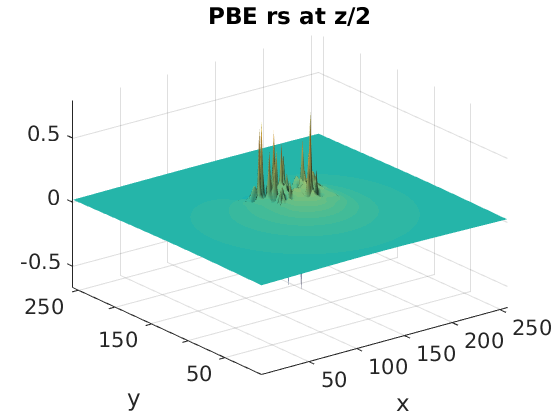}}
 \end{tabular}
 };
 
 \node (Error_PBE_RS_Newt_379_257) at (10,0) {\begin{tabular}{l}
  {\includegraphics[width=0.48\textwidth]{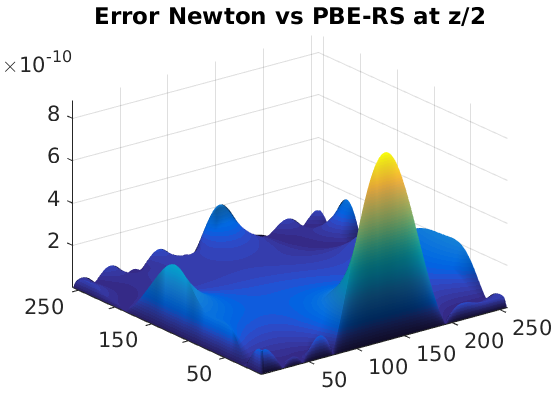}}
 \end{tabular}
 };
 \draw[->,draw=blue,thick] (3.5,0) -- (Error_PBE_RS_Newt_379_257.west);
 
 \end{tikzpicture}}
 \caption{\label{fig:379_Newt_pot_err_257}
 Absolute error between the solutions of the Newton potential sums and the modified Poisson equation for 
 the protein Fasciculin 1.}
\end{figure}

\begin{table}[t]
  \begin{center}%
  \begin{tabular}
  [c]{|c|c|c|c|}%
  \hline
    n         & 97   & 129   & 257   \\
    \hline
  Discrete $L_2$ norm  & $4.1176\times10^{-7}$   & $1.1936\times10^{-7}$ & $3.7003\times10^{-9}$  \\
  \hline
  \end{tabular}
  \caption{\small The discrete $L_2$ norm of the error and the relative error with respect to grid 
  size for the acetazolamide molecule.}
  \label{Tab:Elect_error_Acet}
  \end{center}
\end{table} 
Finally, we consider the protein Fasciculin 1, with 1228 atoms, 
an anti-acetylcholinesterase toxin from the green mamba snake venom \cite{DuMaBoFo:92}.
Again, we compute the electrostatics potential as in the previous case, but with $60\,\textrm{\AA}$ box 
length and a $257^3$ as the minimum Cartesian grid because of a larger molecular size. The results 
provided in the Figure \ref{fig:379_Newt_pot_err_257} and Table \ref{Tab:Elect_error_379} illustrate 
a similar trend of accuracy as in the previous test examples.

\begin{table}[b]
  \begin{center}%
  \begin{tabular}
  [c]{|c|c|c|c|}%
  \hline
    n         & 129   & 193   & 257   \\
    \hline
  Discrete $L_2$ norm  & $1.2919\times10^{-6}$   & $1.7395\times10^{-7}$ & $4.3060\times10^{-8}$  \\
  \hline
  \end{tabular}
  \caption{\small The discrete $L_2$ norm of the error and the relative error with respect to grid 
  size for the protein Fasciculin 1.}
  \label{Tab:Elect_error_379}
  \end{center}
\end{table}

\subsection{Accurate representation of the long-range electrostatic potential by the 
RS tensor}
\label{ssec:Poisson_SR} 

\begin{figure}[t]
\centering
\includegraphics[width=5.1cm]{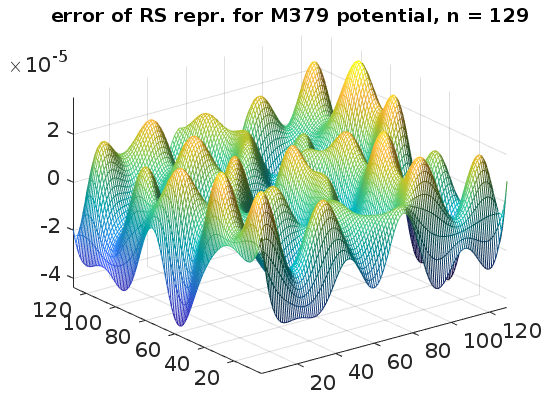} 
\includegraphics[width=5.1cm]{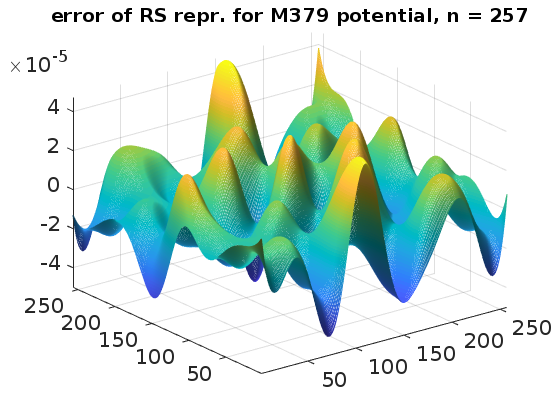} 
\includegraphics[width=5.1cm]{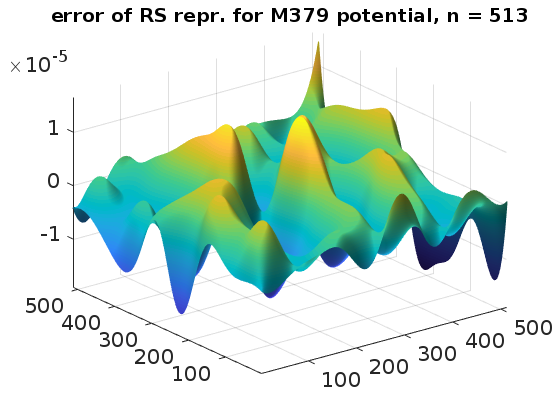}
\caption{\small The error due to the low-rank approximation of the long-range component for the 
379 atomic molecule for $n=129^3$ (left), $n=257^3$ (middle) and $n=513^3$ (right) grids.}
\label{fig:LR_error_06}
\end{figure}

\begin{figure}[t]
\centering
\includegraphics[width=5.8cm]{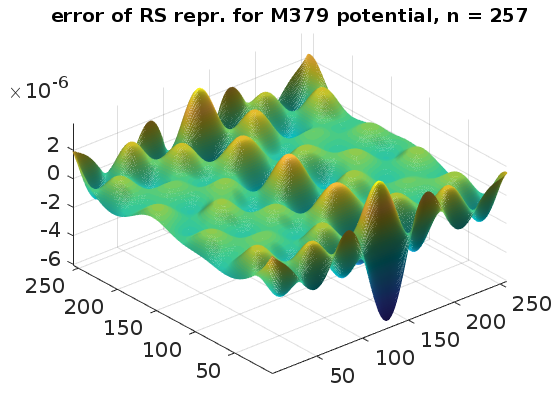}
\caption{\small The error due to the low-rank approximation of the long-range component for the 
379 atomic molecule for $n=257^3$ grid at a lower tolerance.}
\label{fig:LR_error_07}
\end{figure}

Here, we highlight the advantages of the RS tensor format in the low-rank
approximation of the long-range component in the total potential sum. 
For this purpose, the RHOSVD within the multigrid C2T transform \cite{khor-ml-2009} is used 
which provides computation of the low-rank canonical/Tucker tensor representation 
of the long-range part at the asymptotic cost of $\mathcal{O}(Nn)$. 
Here, $N$ is the number of charges in the molecule while $n$ represents the grid dimension in a 
single direction. Figure \ref{fig:LR_error_06} shows an error of $\mathcal{O}(10^{-5})$ of the RS tensor 
format approximation compared with the full size representation at various grid dimensions. 
These data correspond  to the long-range RS-rank equal to $10$, with $\varepsilon$-truncation threshold    
chosen as $\mathcal{O}(10^{-6})$ for the Newton kernel and $\mathcal{O}(10^{-7})$ for the C2T transform. 
Figure \ref{fig:LR_error_07} shows  the error obtained if we take another rank truncation criterion 
$\varepsilon$ of an order less for both the Newton (i.e., $\mathcal{O}(10^{-7})$) 
and for the C2T transform ($\mathcal{O}(10^{-8})$), indicating that
the error will be one order less too, i.e., $\mathcal{O}(10^{-6})$).

\section{Conclusions}\label{sec:Conclusions}

In this paper we demonstrate that the range-separated tensor format is
gainfully applicable for the solution of the PBE for calculation of electrostatics 
in large molecular systems. The efficiency of the new tensor-based regularization scheme for the PBE is 
based on the exceptional properties of the grid-based RS tensor splitting of the Dirac-delta distribution. 
The main computational benefits are due to the  localization of the modified right-hand side within the 
molecular region and automatic maintaining of the continuity of the Cauchy data on the interface.
Another advantage is that our computational scheme only includes solving a single system of algebraic 
equations for the smooth long-range 
(i.e., regularized) part of the collective potential discretized by FDM. The total potential is obtained 
by adding this solution to the directly precomputed low-rank tensor representation for the short-range 
contribution. The various numerical tests illustrate the main properties of the presented scheme. For 
example, it is clear from Figure \ref{fig:Pot_18_sing} that the classical PE model does not accurately 
capture the solution singularities which emanate from the short-range component of the 
total target electrostatic potential in the numerical approximation. In this paper, we emphasize 
that this problem can be efficiently circumvented by applying the range-separated tensor format as a 
solution decomposition technique in order to modify the PBE/PE. In the modified PBE/PE, the Dirac-delta 
distribution is replaced by a smooth long-range function from (\ref{LS_delta}). We thus 
only need to solve for the long-range electrostatic potential numerically and add this solution to 
the short-range component which is computed a priori using the canonical tensor approximation to the 
Newton kernel. The resultant total potential sum is of high accuracy as evident from Figures 
\ref{fig:Newt_pot_err_257} - \ref{fig:379_Newt_pot_err_257} and Tables \ref{Tab:Elect_error_Born} -
\ref{Tab:Elect_error_379}. 
Finally, we summarize that the regularization scheme presented in this paper
has capabilities for various generalizations which can be effectively implemented
with minor changes in the RS tensor decompositions. We notice the following directions:

-- The possibility of PBE computations on much finer grids 
is also an advantage of the proposed approach. Tensor techniques practically may be applied to finer 
grids compared to traditional finite element approaches;

-- The regularization scheme remains verbatim in the case 
of the nonlinear PBE since it requires only the modification of input data for the right-hand side 
in the molecular region $\Omega_m$, where the equation is linear, 
but this does not affect the nonlinearity domain $\Omega_s$; 

-- Our approach allows the efficient calculation of electrostatics under multiple rotations of the 
biomolecule, which is the crucial problem in the numerical modeling of proteins.

  \begin{footnotesize}

  \bibliographystyle{abbrv}
  \bibliography{BSE_Fock_Sums1.bib,Kweyu_refer.bib}

  \end{footnotesize}

\end{document}